\documentclass[10pt]{amsart}
\usepackage[english]{babel}
\usepackage{amssymb}
\usepackage{enumitem}
\usepackage{hyperref}
\usepackage[initials]{amsrefs}
\usepackage[all]{xy}
\SelectTips{cm}{}

\newcommand{\bbQ}{{\mathbb Q}}

\newcommand{\bbZ}{{\mathbb Z}}

\newcommand{\bfG}{{\mathbf G}}
\newcommand{\bfH}{{\mathbf H}}
\newcommand{\bfP}{{\mathbf P}}
\newcommand{\bfQ}{{\mathbf Q}}

\newcommand{\bfA}{{\mathbf A}}

\newcommand{\bfV}{{\mathbf V}}

\newcommand{\calA}{\mathcal{A}}

\newcommand{\calN}{\mathcal{N}}

\newcommand{\calZ}{\mathcal{Z}}

\newcommand{\id}{\operatorname{id}}

\newcommand{\pr}{\operatorname{pr}}

\newcommand{\supp}{\operatorname{supp}}

\newcommand{\Ker}{\operatorname{Ker}}
\newcommand{\Prob}{\operatorname{Prob}}

\newcommand{\Aut}{\operatorname{Aut}}

\newcommand{\PGL}{\operatorname{PGL}}

\newcommand{\Ch}{\operatorname{Ch}}

\newcommand{\overto}[1]{{\buildrel{#1}\over\longrightarrow}}

\newcommand{\acts}{\curvearrowright}
\newcommand{\rk}{\operatorname{rk}}
\newcommand{\gW}[2]{\operatorname{W}_{#1,#2}}
\newcommand{\cW}[1]{\operatorname{Weyl}_{#1}}
\newcommand{\wflip}{w_{\rm flip}}
\newcommand{\wlong}{\operatorname{w}_{\rm long}}
\newcommand{\Sgr}[2]{\mathcal{SG}_{#2}(#1)}
\newcommand{\Quo}[2]{\mathcal{Q}_{#2}(#1)}
\newcommand{\sgr}[1]{\mathcal{SG}(#1)}
\newcommand{\quo}[1]{\mathcal{Q}(#1)}

\newtheorem{theorem}{Theorem}
\newtheorem{lemma}[theorem]{Lemma}

\newtheorem{proposition}[theorem]{Proposition}

\theoremstyle{definition}
\newtheorem{definition}[theorem]{Definition}

\newtheorem{examples}[theorem]{Examples}

\begin{document}

\title{Boundaries, Weyl groups, and Superrigidity}

\author{Uri Bader}
\address{Technion, Haifa}
\email{uri.bader@gmail.com}

\author{Alex Furman}
\address{University of Illinois at Chicago, Chicago}
\email{furman@math.uic.edu}

\thanks{U.B. and A.F. were supported in part by the BSF grant 2008267.}
\thanks{U.B was supported in part by the ISF grant 704/08.}
\thanks{A.F. was supported in part by the NSF grants DMS 0905977.}


\maketitle

\begin{abstract}
This note describes a unified approach to several superrigidity results, old and new,
concerning representations of lattices into simple algebraic groups over local fields.
For an arbitrary group $\Gamma$ and a boundary action $\Gamma\acts B$
we associate certain generalized Weyl group $\gW{\Gamma}{B}$ and show that any
representation with a Zariski dense unbounded image in a simple algebraic group,
$\rho:\Gamma\to \bfH$,
defines a special homomorphism $\gW{\Gamma}{B}\to\cW{\bfH}$.
This general fact allows to deduce the aforementioned superrigidity results.
\end{abstract}

\bigskip

\subsection*{Introduction} 
\label{sub:introduction}\hfill{}\\
This note describes some aspects of a unified approach to a family of
"higher rank superrigidity" results, based on a notion of a generalized Weyl group.
While this approach applies equally well to representations of lattices
(as in the original work of Margulis \cite{Margulis-ICM}),
and to measurable cocycles (as in the later work of Zimmer \cite{Zimmer-csr}),
in this note we shall focus on representations only.
Yet, it should be emphasized that our techniques do not involve any cocompactness,
or integrability assumptions on lattices, and their generalizations to general measurable
cocycles are rather straightforward.
Hereafter we consider representations into simple algebraic groups;
some other possible target groups are discussed in \cite{Bader+Furman:hyp}, \cite{Bader+Furman+Shaker}.

Let $k$ be a local field, and $\bfH$ denote the locally compact group of $k$-points of some
connected adjoint $k$-simple $k$-algebraic group.
Consider representations $\rho:\Gamma\to \bfH$ with Zariski dense and unbounded image,
where $\Gamma$ is some discrete countable group.
We shall outline a unified argument showing that for the following groups $G$
all lattices $\Gamma$ in $G$ have the property that such a representation $\rho:\Gamma\to\bfH$
can occur only as a restriction of a continuous homomorphism $\bar\rho:G\to \bfH$;
this includes:
\begin{itemize}
	\item[(a)]
	$G=\bfG$ is the group of $\ell$-points of a connected $\ell$-simple $\ell$-algebraic group of
	$\rk_\ell(G)\ge 2$ where $\ell$ is a local field (Margulis \cite{Margulis-ICM}, \cite[\S VII]{Margulis-book}),
	\item[(b)]
	$G=G_1\times G_2$ for $G_1,G_2$ general locally compact groups,
	where $\Gamma$ is assumed to be irreducible (cf. \cite{Monod-CAT0}, \cite{GKM}, \cite{BMZ}),
	\item[(c)]
	$G=\Aut(X)$ where $X$ is an $\widetilde{A}_2$-building and $G$ has finitely many orbits for its action
    on the space of chambers of $X$, $\Ch(X)$.
\end{itemize}
New implications of these results include non-linearity of the exotic $\widetilde{A}_2$-groups
(deduced from (c)), and arithmeticity vs. non-linearity dichotomy for irreducible lattices in
products of topologically simple groups, as in \cite{Monod-arith},
but with integrability assumptions removed.
Cocycle versions of the above results cover more new ground.

\subsection*{I. Boundaries and Weyl groups} 
\label{sub:boundaries_and_weyl_map}\hfill{}\\
We start from some constructions related to the source group $G$.
In rigidity theory, a boundary of a lcsc group $G$ is an auxiliary measure space $(B,\nu)$
equipped with a measurable, measure-class preserving action of $G$,
satisfying additional conditions that imply existence of
measurable $\Gamma$-equivariant maps from $B$ to some compact homogeneous $\bfH$-spaces.
We shall work with the following:
\begin{definition}
A measurable $G$-space $(B,\nu)$ is a $G$-\emph{boundary} if:
\begin{itemize}
	\item[(B1)] the action $G\acts B$ is \emph{amenable} in the sense of Zimmer \cite{Zimmer-amen},
	\item[(B2)] the projection $\pr_1:(B\times B,\nu\times\nu)\to (B,\nu)$ is
    \emph{ergodic with Polish coefficients} (for short, EPC), as defined below.
\end{itemize}
\end{definition}

We say that the action $G\acts (X,\mu)$ is \emph{ergodic with Polish coefficients} (\emph{EPC})
if for every isometric $G$-action on any Polish metric space $(U,d)$, every measurable $G$-equivariant map $F:X\to U$
is $\mu$-a.e. constant. 
We say that a measure class preserving $G$-\emph{map} $\pi:(X,\mu)\to(Y,\nu)$ is \emph{EPC}
(or that $X$ is \emph{EPC relatively to $Y$}, when the map $\pi$ is understood)
if for every $G$-action by fiber-wise isometries on a measurable field
of Polish metric spaces $\{ (U_v,d_v) \}_{v\in V}$ over a standard Borel $G$-space $V$,
every measurable $G$-map $F:X\to U=\bigsqcup_{v\in V} U_v$ descends to $Y$, i.e.
there exists a measurable $f:Y\to U$ so that $\mu$-a.e. $F$ coincides
with $f\circ \pi$:
\[
	\xymatrix{ (X,\mu) \ar[r]^{F} \ar[d]_{\pi} & U \ar[d] \\
	           (Y,\nu) \ar[r] \ar@{.>}[ur]^{f} & V
	}
\]
The relative EPC property (B2) implies the (absolute) EPC property,
which in turn implies ergodicity with unitary coefficients for $G\acts B\times B$ and $G\acts B$.
Hence $G$-boundaries are also strong boundaries in the sense of Burger-Monod \cite{Burger+Monod}.
Moreover, it can be shown that:
\begin{enumerate}
	\item
	If $\eta$ is a symmetric spread out generating probability measure on a lcsc group $G$, then the Poisson-Furstenberg
	boundary $(B,\nu)$ of $(G,\eta)$ is a $G$-boundary (further strengthening \cite{Kaimanovich}).
	\item
	If $\bfG$ is a simple algebraic $\ell$-group and $\bfP<\bfG$ is a minimal parabolic,
	then $B=\bfG/\bfP$ with the Haar measure class is a $\bfG$-boundary.
	\item
	If $(B_i,\nu_i)$ are $G_i$-boundaries, for $i=1,2$,
    then $(B,\nu)=(B_1,\nu_1)\times (B_2,\nu_2)$
	is a $G$-boundary for the product $G=G_1\times G_2$.
	\item
	If $\Gamma<G$ is a lattice, then every $G$-boundary is also a $\Gamma$-boundary.
\end{enumerate}
\noindent{Given} a $G$-boundary $(B,\nu)$ consider the group
$\Aut_G(B\times B)$ of all measure class preserving automorphisms of $(B\times B,\nu\times\nu)$
which are equivariant with respect to the diagonal $G$-action.
The flip involution $\wflip:(b,b')\mapsto (b',b)$ is an obvious example
of such a map.
A \emph{generalized Weyl group} $\gW{G}{B}$ associated to a choice of a $G$-boundary
is $\Aut_G(B\times B)$, or a subgroup of $\Aut_G(B\times B)$ containing the flip $\wflip$.
If $\Gamma$ is a lattice in $G$, any $G$-boundary $B$ is also a $\Gamma$-boundary, and we can take
\[
	\gW{\Gamma}{B}=\Aut_G(B\times B).
\]
We view $B$ and $\gW{G}{B}$ as auxiliary objects, associated (not in a unique way) to $G$,
and encoding its implicit symmetries.
Non-amenable groups have non-trivial boundaries, so their generalized Weyl groups always contain
$\{\id, \wflip\}\cong \bbZ/2\bbZ$.
Presence of additional elements can be viewed as an indication of "higher rank phenomena",
as in the following examples:
\begin{examples}\label{E:Weylgroups}
\begin{itemize}
	\item[(a)]
	Let $\mathbf{G}$ be a non-compact simple algebraic group and $B=\bfG/\bfP$ its flag variety.
	Then $B\times B\cong \bfG/\calZ_\bfG(\bfA)$ as measurable $G$-spaces, and the
	generalized Weyl group coincides with the classical one:
	\[
		\gW{\bfG}{B}= \Aut_\bfG(\bfG/\calZ_\bfG(\bfA))
		=\calN_\bfG(\bfA)/\calZ_\bfG(\bfA)=\cW{\bfG}.
	\]
	In particular, $\cW{\PGL_n(\ell)}\cong S_n$.
	Note that $\cW{\bfG}\not\cong \bbZ/2\bbZ$ iff $\rk_\ell(\bfG)\ge 2$.
	\item[(b)]
	Let $G=G_1\times G_2$ be a product of non-amenable factors,
	take $B=B_1\times B_2$ and
	$\gW{G}{B}=\{\id,w_1,w_2,\wflip\}$ where the elements $w_i$ are the flips of the $B_i$-coordinates.
    Note that $\gW{G}{B}\cong (\bbZ/2\bbZ)^2\le \gW{G_1}{B_1}\times\gW{G_2}{B_2}$.
	\item[(c)]
	Let $G=\Aut(X)$ be an $\widetilde{A}_2$-group, the Poisson-Furstenberg boundary for the simple random walk
	can be realized on the space of chambers of the associated spherical building, $B=\Ch(\partial X)$,
	with $\gW{G}{B}\cong S_3$.
	This generalizes the classical case $\cW{\PGL_3(\bbQ_p)}\cong S_3$.
\end{itemize}
\end{examples}


\subsection*{II. The homomorphism between Weyl groups} 
\label{sub:the_bfh_gate_of_an_ergodic_gamma_space}\hfill{}\\
Next we turn to the target group $\bfH$. 
Denote by $\bfA<\bfP<\bfH$ the ($k$-points of) a maximal $k$-split torus and minimal parabolic subgroup containing it.
The diagonal $\bfH$-action on $\bfH/\bfP\times\bfH/\bfP$ has finitely many orbits, indexed by $\cW{\bfH}$
(Bruhat decomposition) with a unique full-dimensional orbit,
corresponding to the long element $\wlong\in \cW{\bfH}$:
\begin{equation}\label{e:big-cell}
	\bfH/\calZ_\bfH(\bfA)\subset \bfH/\bfP\times\bfH/\bfP,\qquad h\calZ_\bfH(\bfA)\mapsto (h\bfP,h\wlong\bfP).
\end{equation}
In the case of $\bfH=\PGL_n(k)$ the groups $\calZ_\bfH(\bfA)<\bfP$ correspond to the diagonal and the upper triangular subgroups,
$\bfH/\bfP$ is the space of flags $(E_1\subset E_2\subset\dots\subset E_n)$ where $E_j$ is a $j$-dimensional subspace of $k^n$,
$\bfH/\calZ_\bfH(\bfA)$ is the space of $n$-tuples of one-dimensional subspaces $(\ell_1,\dots,\ell_n)$
with $\ell_1\oplus\cdots\oplus \ell_n=k^n$, and (\ref{e:big-cell}) is
\[
	(\ell_1,\dots,\ell_n)\mapsto ((\ell_1,\ell_1\oplus\ell_2,\dots),
	(\ell_n,\ell_n\oplus\ell_{n-1},\dots)).
\]
Here $\cW{\PGL_n(k)}=S_n$ acts by permutations on $(\ell_1,\dots,\ell_n)$ with the long element
$\wlong:(\ell_1,\dots,\ell_n)\mapsto(\ell_n,\dots,\ell_1)$.

We have the following general result:
\begin{theorem}\label{T:gates+Weyl}
	Let $\Gamma$ be a countable group, $\rho:\Gamma\to\bfH$ a homomorphism with unbounded and Zariski dense image,
	and $(B,\nu)$ be a $\Gamma$-boundary. Then
	\begin{itemize}
		\item[{\rm (i)}]
		There exists a unique measurable $\Gamma$-equivariant map $\phi:B\to \bfH/\mathbf{P}$,
		\item[{\rm (ii)}]
		The map $\phi\times\phi:B\times B\to \bfH/\mathbf{P}\times \bfH/\mathbf{P}$
		factors through the embedding
		\[
			\bfH/\mathcal{Z}_\bfH(\mathbf{A})\to \bfH/\mathbf{P}\times \bfH/\mathbf{P},
		\]
		\item[{\rm (iii)}]
		There exists a homomorphism $\pi:\gW{\Gamma}{B}\to \cW{\bfH}$ with
		\[
			\qquad\pi(\wflip)=\wlong,\qquad
			(\phi\times\phi)\circ w=\pi(w)\circ (\phi\times\phi)\qquad
			(w\in \gW{\Gamma}{B})
		\]
		as measurable maps $B\times B\to \bfH/\mathcal{Z}_\bfH(\mathbf{A})\subset \bfH/\mathbf{P}\times \bfH/\mathbf{P}$.
	\end{itemize}
\end{theorem}
%
Let us sketch the main ingredients in the proof of this result.
Consider $k$-algebraic actions of $\bfH$ on ($k$-points of) $k$-algebraic varieties, or $\bfH$-\emph{varieties}
for short.
Orbits of such actions are locally closed, and the actions are \emph{smooth}, in the sense
that the space of orbits is standard Borel.
Since orbits of ergodic actions cannot be separated it follows (\cite[2.1.10]{Zimmer-book})
that given an ergodic measure class preserving action $\Gamma\acts (S,\sigma)$
any measurable $\Gamma$-equivariant map $\phi:S\to \bfV$ into a $\bfH$-variety,
takes values in a single $\bfH$-orbit $\bfH v\subset \bfV$.
This can be used to show the following general
\begin{lemma}\label{L:gate}
	Let $\Gamma\acts (S,\sigma)$ be an ergodic measure-class preserving action, and $\rho:\Gamma\to \bfH$ a homomorphism.
	There exists $\phi_0:S\to \bfV_0\cong \bfH/\bfH_0$, where $\bfH_0<\bfH$ is an algebraic subgroup, so that
	for any measurable $\Gamma$-map  $\phi:S\to \bfV$ to a $\bfH$-variety, there is a unique $\bfH$-map
	$f:\bfH/\bfH_0\to\bfV$ so that $\phi=f\circ \phi_0$ a.e. on $S$.
\end{lemma}
In other words, $\phi_0:S\to \bfH/\bfH_0$ is an initial object in the category consisting of measurable
$\Gamma$-maps $\phi:S\to\bfV$ to $\bfH$-varieties, where morphisms between $\phi_i:S\to\bfV_i$ ($i=1,2$)
are algebraic $\bfH$-maps $f:\bfV_1\to\bfV_2$ with $\phi_2=f\circ \phi_1$.
The initial object described in Lemma~\ref{L:gate} is necessarily unique up to $\bfH$-automorphisms
of $\bfH/\bfH_0$, i.e. up to the action of $\calN_{\bfH}(\bfH_0)/\bfH_0$ from the right.
This universality gives a homomorphism from the group of measure-class preserving automorphisms
of the $\Gamma$-action to the latter group
\[
	\Aut_\Gamma(S,[\sigma])\to \Aut_\bfH(\bfH/\bfH_0)=\calN_{\bfH}(\bfH_0)/\bfH_0.
\]
We apply this construction to the diagonal $\Gamma$-action on $S=B\times B$, and use it
to deduce Theorem~\ref{T:gates+Weyl} from the following
\begin{theorem}
Let $\rho:\Gamma\to \bfH$ be as above, and $(B,\nu)$ be a $\Gamma$-boundary. Then:
\begin{itemize}
	\item[{\rm (i)}]
	The initial $\bfH$-object for $\Gamma\acts (B,\nu)$ is $\phi:B\to \bfH/\bfP$.
	\item[{\rm (ii)}]
	The initial $\bfH$-object for the diagonal action $\Gamma\acts (B\times B,\nu\times \nu)$ is
	\[
		\phi\times\phi:B\times B\to \bfH/\calZ_\bfH(\calA)\subset \bfH/\bfP\times\bfH/\bfP.
	\]
\end{itemize}
\end{theorem}
Note that condition (B1) of amenability of $\Gamma\acts B$ yields a measurable $\Gamma$-map $\Phi:B\to\Prob(\bfH/\bfP)$.
Claim (i) asserts that $\Phi(b)=\delta_{\phi(b)}$ are Dirac measures, and that $B$ admits no $\Gamma$-maps
to $\bfH/\bfH_0$ where $\bfH_0$ is a proper algebraic subgroup of $\bfP$.
The proof of (ii) relies on the relative EPC property of $B\times B\to B$ in showing
that $B\times B$ has no $\Gamma$-maps to $\bfH/\bfH_0$ with $\bfH_0$ a proper subgroup of $\calZ_\bfH(\bfA)$.


\subsection*{III. Galois correspondence} 
\label{sub:the_galois_correspondence}\hfill{}\\
Let $B$ be a set and $W$ be a group acting on $B\times B$.
This very general datum alone defines an interesting structure that we shall now briefly describe
(see \cite{Bader+Furman:hyp} for more details).

Consider possible quotients $p:B\to p(B)$ of $B$, or rather
their equivalence classes determined by the pull-back of the Boolean algebra
from $p(B)$ to $B$.
Denote by $\quo{B}$ the collection of all such (classes of) quotients, ordered by
$p_1\le p_2$ if $p_1=j\circ p_2$ for some $j:p_2(B)\to p_1(B)$ or, equivalently, by
inclusion of the corresponding Boolean algebras.
Let $\sgr{W}$ denote all subgroups of $W$, ordered by inclusion,
and consider maps $p\mapsto W_p$ and $V\mapsto p^V$
between $\quo{B}$ and $\sgr{W}$ defined as follows.
Given $p\in \quo{B}$ consider the map $p_1:B\times B\overto{\pr_1} B\overto{p} p(B)$
and define
\[
	W_p=\{ w\in W \mid p_1\circ w=p_1 \}.
\]
Given a subgroup $V\le W$ define the quotient $p^V:B\to p^V(B)$ to be the finest one with $V\le W_p$.
Then the maps $p\mapsto W_p$, $V\mapsto p^V$, between $\quo{B}$ and $\sgr{W}$,
viewed as partially ordered sets, are order-reversing and satisfy $V\le W_p$ iff $p\le p^V$.
A pair of order-reversing maps between posets with above property,
forms an abstract \emph{Galois correspondence}; one of the formal consequences
of such a setting is that one can define the following operations of taking a \emph{closure}
in $\quo{B}$, $\sgr{W}$:
\[
	\overline{V}:= W_{(p^V)},\quad \overline{p}:=p^{(W_p)}\qquad\textrm{satisfying}\qquad
	V\le \overline{V}=\overline{\overline{V}},\quad p\le \overline{p}=\overline{\overline{p}}.
\]
It follows that the collections of closed objects in $\quo{B}$ and $\sgr{W}$
\[
	\Quo{B}{W}=\{\overline{p}:B\to\overline{p}(B)\}, \qquad \Sgr{W}{B}=\{\overline{V} \mid V\le W\}
\]
form sub-lattices\footnote{ The term \emph{lattice} here refers to a partially ordered set,
where any two elements $x,y$ have a \emph{join} $x\vee y$ and \emph{meet} $x\wedge y$.}
of $\quo{B}$, $\sgr{W}$, on which the above Galois correspondence is an order-reversing isomorphism.

These constructions can be carried over to the measurable setting,
where $(B,\mathcal{B},\nu)$ is a measure space, $\quo{B}$ consists of measurable quotients
(equivalently, complete sub-$\sigma$-algebras of $\mathcal{B}$),
and $W$ is assumed to preserve the measure class $[\nu\times\nu]$.
\begin{examples}\label{E:Galois}
	With $W=\gW{G}{B}$ acting on a double of a $G$-boundary $(B,\nu)$,
	as in Examples~\ref{E:Weylgroups} we get:
\begin{itemize}
	\item[(a)]
	Let $\bfG$ be a simple algebraic group, $B=\bfG/\bfP$ and $W=\cW{\mathbf{G}}$. Then
	\[
		\Quo{\bfG/\bfP}{\cW{\bfG}}=\{\bfG/\bfP\to \bfG/\bfQ \mid \bfP\le \bfQ\le \bfG \textrm{\ parabolic}\}
	\]
	corresponding to Weyl groups $\cW{{\rm Levi}(\bfQ)}$ embedded in $\cW{\bfG}$.
	The lattice is a $\rk_\ell(G)$-dimensional cube.
    \item[(b)]
	Let $G=G_1\times G_2$, $B=B_1\times B_2$ and
	$W=\{\id, w_{1},w_{2}, \wflip\}$.
	Then the non-trivial closed objects are the two factors $B_i$, corresponding to $\{\id, w_{i}\}$.
	\item[(c)]
	For an $\widetilde{A}_2$-group $G=\Aut(X)$, the non-trivial closed quotients of
	$B=\Ch(\partial X)$ are face maps, corresponding to $\{ \id, (1,2)\}$, $\{\id, (2,3)\}$ in $S_3$.
\end{itemize}
\end{examples}
Let $\Gamma<G$ be lattice where $G$ is one of the above examples,
view the $G$-boundary $(B,\nu)$ as a $\Gamma$-boundary and take $\gW{\Gamma}{B}=W$.
The Galois correspondence above was determined by $W\acts B\times B$ alone (without any reference to a $G$-action);
so the concepts of closed subgroups and closed quotients remain unchanged.

Given an unbounded Zariski dense representation $\rho:\Gamma\to\bfH$, we consider
the associated map $\phi:B\to \bfH/\bfP$ and homomorphism  $\pi:\gW{\Gamma}{B}\to \cW{\bfH}$
as in Theorem~\ref{T:gates+Weyl}.
Then $\Ker(\pi)$ is a normal subgroup in $\gW{\Gamma}{B}$ which is also closed in the above sense,
and $\phi$ factors through a closed quotient corresponding to $\Ker(\pi)$.
\begin{proposition}\label{P:semi-final}
	Let $\rho:\Gamma\to\bfH$, $\phi:B\to \bfH/\bfP$, $\pi:\gW{\Gamma}{B}\to\cW{\bfH}$ be as above:
	\begin{itemize}
		\item[{\rm (a)}]
		If $\Gamma<\bfG$ is a lattice in a simple algebraic group,
		then $\pi:\gW{\Gamma}{B}=\cW{\bfG}\to\cW{\bfH}$ is injective and $\pi(\wlong^{(\bfG)})=\wlong^{(\bfH)}$.
		\item[{\rm (b)}]
		If $\Gamma<G=G_1\times G_2$, $B=B_1\times B_2$, and
		$\pi:\gW{\Gamma}{B}=(\bbZ/2\bbZ)^2\to\cW{\bfH}$ is non-injective,
		then $\phi$ factors through some
		\[
			\phi:B\ \overto{\pr_i}\ B_i\ \overto{\phi_i}\ \bfH/\bfP.
		\]
		where $\phi_i:B_i\to \bfH/\bfP$ is a measurable $\Gamma$-map for some $i=1,2$.
		\item[{\rm (c)}]
		If $\Gamma<G=\Aut(X)$ is an $\widetilde{A}_2$-group, then $\pi:S_3\to\cW{\bfH}$ is injective
		with $\pi((1,3))=\wlong^{(\bfH)}$.
	\end{itemize}
\end{proposition}
Case (b) follows from the classification of closed subgroups and corresponding quotients for products;
while (a) and (c) are consequences of a general fact (\cite{Bader+Furman:hyp}) that irreducible Coxeter groups,
such as $\cW{\bfG}$, have no non-trivial normal \emph{special subgroups} (same as closed subgroups in our context).


\subsection*{IV. The final step} 
\label{sub:the_end_game}\hfill{}\\
Proposition~\ref{P:semi-final} already suffices to deduce some superrigidity
results of the type \emph{``certain $\Gamma$ admits no unbounded Zariski dense homomorphisms to certain $\bfH$''}.
For example, this is the case if $\bfH$ is a simple $k$-algebraic group with $\rk_k(\bfH)=1$,
while $\Gamma$ is an exotic $\widetilde{A}_2$-group or a lattice in a simple $\ell$-algebraic group $\bfG$ with $\rk_\ell(\bfG)\ge 2$.
Embedding of $\cW{\bfG}$ in $\cW{\bfH}$ preserving the long element, can be ruled out in many other
cases, such as $\bfG=\PGL_4(\ell)$ and $\bfH=\PGL_3(k)$.

However, in the context of Proposition~\ref{P:semi-final} one has more precise information:
\begin{itemize}
	\item[{\rm (a)}]
	If $\Gamma$ is a lattice in a simple algebraic group $\bfG$,
	then $\pi:\gW{\Gamma}{B}=\cW{\bfG}\to\cW{\bfH}$ is an isomorphism of Coxeter groups.
	\item[{\rm (b)}]
		If $\Gamma<G=G_1 \times G_2$, $B=B_1\times B_2$, then
		$\pi:\gW{\Gamma}{B}=(\bbZ/2\bbZ)^2\to\cW{\bfH}$ has $\operatorname{Im}(\pi)\cong \bbZ/2\bbZ$.
	\item[{\rm (c)}]
		If $\Gamma<G=\Aut(X)$ is an $\widetilde{A}_2$-group, then $\pi:S_3\to\cW{\bfH}$ is an
		isomorphism of Coxeter groups.
	\end{itemize}
The proof of these claims relies on classification of measurable $\Gamma$-factors of $\Gamma$-boundaries
proved by Margulis \cite{Margulis-Factor}, Bader-Shalom \cite{Bader+Shalom}, and Shalom-Steger \cite{Shalom+Steger}, respectively.
Let us turn to consequences of such statements.
\begin{theorem}
	Let $\Gamma<G=G_1\times G_2$ be a lattice with $\pr_i(\Gamma)$ dense in $G_i$ ($i=1,2$), and
	$\rho:\Gamma\to\bfH$ a Zariski dense unbounded representation.
	Then $\rho$ extends to $G$ and factors through a continuous homomorphism $\rho_i:G_i\to\bfH$
	of a factor.
\end{theorem}
This follows from the fact that $\pi:\gW{\Gamma}{B}\to\cW{\bfH}$ has two element image,
Proposition~\ref{P:semi-final}.(b), and the following general
\begin{lemma}
	Let $G_1\acts (B_1,\nu_1)$ be a measure class preserving action of some locally compact group $G_1$,
	$\Gamma$ any group, $p:\Gamma\to G_1$ a homomorphism with dense image; and let
	$\phi_1:B_1\to \bfH/\bfP$ a measurable $\Gamma$-equivariant map.
	Then there is a continuous homomorphism $\bar\rho_1:G_1\to \bfH$ so that $\rho=\bar\rho_1\circ p$.
\end{lemma}
The proof of this lemma utilizes the enveloping semigroup of $\bfH\acts\bfH/\bfP$, which can be 
identified as the \emph{quasi-projective transformations} of $\bfH/\bfP$, introduced in \cite{Furstenberg}.

The final treatment of case (a) (the main case of Margulis's superrigidity) and the non-linearity
result for exotic $\widetilde{A}_2$-groups (as in (c)) are deduced by a reduction
to some results of Tits on buildings.
If $G=\bfG$ is a simple $\ell$-algebraic group let $\Delta=\Delta_\bfG$ denote
the spherical building of $\bfG$, if $G=\Aut(X)$ is an $\widetilde{A}_2$-group,
let $\Delta=\partial X$ denote the spherical building associated to the Affine building $X$.
Let $\Delta'=\Delta_\bfH$ denote the spherical building of $\bfH$.
For a building $\Delta$ denote by $\Ch(\Delta)^{(2)}$ the subspace of $\Ch(\Delta)\times\Ch(\Delta)$
consisting of pairs of opposite chambers.

\begin{theorem}
	Let $\Delta$ and $\Delta'$ be two spherical buildings of rank $\ge 2$
	with the same Weyl group.
	Let $\nu$ be a probability measure on $\Ch(\Delta)$ with $\supp(\nu)=\Ch(\Delta)$
	and $(\nu\times\nu)(\Ch(\Delta)^{(2)})=1$,
	and let $\phi:\Ch(\Delta)\to\Ch(\Delta')$ be a measurable map so that
	$(\phi_*\nu\times\phi_*\nu)(\Ch(\Delta')^{(2)})=1$ and
	\[
		(\phi\times\phi)\circ w=w\circ (\phi\times\phi)\qquad (w\in W_{\Delta}).
	\]
	Then there exists an imbedding of buildings $\Delta\to \Delta'$ which induces a continuous map
    $\Ch(\Delta)\to\Ch(\Delta')$ which agrees $\nu$-a.e. with $\phi$.
\end{theorem}

The proof is now completed by invoking results of Tits, on reconstructing $\bfG$ from $\Delta_\bfG$,
and $X$ from $\partial X$.

%



\end{document}